\documentclass{amsart}
\usepackage{pstricks,pst-node,amssymb}

\title{Quantum Bruhat graph and Schubert polynomials}

\author{Alexander Postnikov}
\date{June~7, 2002}
\address{Department of Mathematics, M.I.T., Cambridge, MA 02139}
\email{apost@math.mit.edu}
\urladdr{http://www.math-mit.edu/\mytilde apost/}

\thanks{The author was supported in part by NSF grant DMS-0201494.}

\newtheorem{theorem}{Theorem}

\newtheorem{corollary}[theorem]{Corollary}

\newtheorem{lemma}[theorem]{Lemma}

\def\Z{\mathbb{Z}}
\def\C{\mathbb{C}}
\def\Q{\mathbb{Q}}

\def\mytilde{\kern-.015in\hbox{\lower.03in\hbox{\~{}}}\kern-.01in}
\def\<{\left<}
\def\>{\right>}
\def\QH{\mathrm{QH}}
\def\HH{\mathrm{H}}

\def\wnot{{w_\mathrm{o}}}
\def\T{\mathcal{T}}
\def\H{\mathcal{H}}
\def\E{\mathcal{E}}
\def\S{\mathfrak{S}}
\def\Fl{\mathrm{Fl}}
\def\s{\sigma}
\def\ds{\displaystyle}
\def\weight{\mathrm{weight}}
\def\dmin{d_\mathrm{min}}
\def\h{\mathfrak{h}}
\def\End{\mathrm{End}}

\psset{unit=1pt, arrowsize=2pt 1, linewidth=0.5pt, linecolor=blue}

\begin{document}

\begin{abstract}
The quantum Bruhat graph, which is an extension of the graph formed
by covering relations in the Bruhat order, is naturally related to the quantum
cohomology ring of $G/B$.  We enhance a result of Fulton and Woodward
by showing that the minimal monomial in the quantum parameters that
occurs in the quantum product of two Schubert
classes has a simple interpretation in terms of directed paths in this graph.

We define path Schubert polynomials, which are quantum cohomology
analogues of skew Schubert polynomials recently introduced by 
Lenart and Sottile.  They are given by sums over paths in the quantum 
Bruhat graph of type $A$.  The 3-point Gromov-Witten invariants 
for the flag manifold are expressed in terms of these polynomials.  This 
construction gives a combinatorial description for the set of all 
monomials in the quantum parameters that occur in the quantum product 
of two Schubert classes.
\end{abstract}

\maketitle

\section{Introduction}
\label{sec:intro}

The famous Littlewood-Richardson coefficients are the structure constants
of the cohomology ring of the Grassmannian in the basis of Schubert classes.
It is an open question to give a combinatorial interpretation of the 
generalized Littlewood-Richardson coefficients for the flag manifold and, 
more generally, for the homogeneous space $G/P$.
This article discusses the 3-point {\it Gromov-Witten invariants}, 
which are the structure constants of the small {\it quantum cohomology ring}.
They extend the Littlewood-Richardson coefficients in an orthogonal 
``quantum'' direction. 

Fulton and Woodward~\cite{FW} described the set of minimal 
monomials $q^d$ in the quantum parameters that occur in the quantum 
product $\s_u*\s_u$ of two Schubert classes in $\QH^*(G/P)$.  
We demonstrate that the problem of multiplying the Schubert
classes in the quantum cohomology ring is naturally related to the study
of paths in the {\it quantum Bruhat graph\/} from~\cite{BFP}.
It is a directed graph is obtained by adding some ``quantum'' 
edges to the graph formed by covering relations in the Bruhat order.
The quantum Bruhat graph encodes the terms that appear in the quantum 
Chevalley-Monk formula.
We show, for $G/B$, that there is a {\it unique\/} minimal monomial $q^d$ 
that occurs in the quantum product $\s_u*\s_v$.
It is equal to the weight of any {\it shortest\/} directed path from 
$u$ to $\wnot v$ in the quantum Bruhat graph.  
All such shortest paths should have the same weight.

Lenart and Sottile~\cite{LS} recently defined {\it skew Schubert polynomials}.
The expansion coefficients of these polynomials in the basis of usual Schubert 
polynomials are equal to the generalized Littlewood-Richardson coefficients
for the flag manifold.  On the other hand, in~\cite{P3}, we
introduced toric Schur polynomials, whose expansion coefficients
in the basis of usual Schur polynomials are the Gromov-Witten invariants 
for the Grassmannian.

In this article, we put these two approaches together.
We define the {\it path Schubert polynomials}, whose expansion coefficients
in the basis of usual Schubert polynomials give the Gromov-Witten invariants 
for the flag manifold.
The path Schubert polynomials can be defined as sums over directed paths
in the quantum Bruhat graph on the symmetric group. 
The proof of the result on their Schubert-expansion readily follows 
from a combination of the Cauchy identity,
the elementary quantization rule from~\cite{FGP},
and the quantum Pieri formula.

This construction implies a combinatorial description for 
the set of {\it all\/} monomials $q^d$ in the quantum parameters that occur,
with non-zero coefficients, in the quantum product $\s_u*\s_v$ of two 
Schubert classes in the quantum cohomology of the flag manifold.  
A monomial $q^d$ occurs in $\s_u*\s_v$ if and only if there is a directed 
path of weight $q^d$ from the vertex $u$ to the 
the vertex $\wnot v$ in the quantum Bruhat graph that 
satisfy certain condition.  This solves a problem posed by Fulton 
and Woodward~\cite{FW}, for the flag manifold.

\section{Quantum cohomology of $G/B$}

In this section we discuss the quantum cohomology ring
of the generalized flag manifold $G/B$.  
First, we briefly recall some general Lie-theoretic notation, 
see~\cite{Hum} for more details.  Then we remind a few facts related to 
Schubert classes and quantum cohomology, and formulate the quantum 
Chevalley-Monk formula, see~\cite{FW}.
This formula leads to the quantum Bruhat graph that was introduced 
in~\cite{BFP}.  We show that the minimal monomial $q^d$ that occurs 
in the quantum product of two Schubert classes has a simple 
interpretation in terms of directed paths in the quantum Bruhat graph.   
\medskip

Let $G$ be a simply connected complex semi-simple Lie group together 
with a Borel subgroup $B$.
The homogeneous space $G/B$ is a compact complex manifold.

Let $\h$ be the {\it Cartan subalgebra\/} in the Lie algebra $\mathfrak{g}$ 
of $G$, and let $\Phi\subset\h^*$ be the {\it root system\/} associated 
with $G$.  For a root $\alpha\in\Phi$, let $h_\alpha\in \h$ 
denote the corresponding {\it coroot\/} in the dual root system 
$\Phi^\vee$.  The {\it Weyl group\/} $W$ is generated by the reflections 
$s_\alpha\in\End(\h^*)$, $\alpha\in\Phi$, given by
$s_\alpha:\lambda\mapsto \lambda-\lambda(h_\alpha)\, \alpha$.
If $\<\alpha\mid\beta\>$ is a $W$-invariant inner product on $\h^*$, 
then $\h$ can be identified with $\h^*$ via this inner product
and $h_\alpha=2\alpha/\<\alpha\mid\alpha\>$.
The choice of Borel subgroup $B$ determines the set of {\it positive roots\/}
$\Phi_+\subset \Phi$, the basis of {\it simple roots\/} 
$\alpha_1,\dots,\alpha_r$ in the root system~$\Phi$,
and the basis of simple coroots $h_{\alpha_1},\dots,h_{\alpha_r}$ 
in the dual root system $\Phi^\vee$. 
The {\it fundamental weights\/} $\lambda_1,\dots,\lambda_r\in\h^*$ are
defined as elements of the basis dual to 
$\{h_{\alpha_1},\dots,h_{\alpha_r}\}$, 
i.e., $\lambda_i(h_{\alpha_j})=\delta_{ij}$.
The Weyl group $W$ is generated by the set of 
{\it Coxeter generators\/} $s_i=s_{\alpha_i}$.
The {\it length\/} $\ell(w)$ of an element $w\in W$ is the 
minimal number of generators $s_i$ in a decomposition for~$w$.

The {\it cohomology ring\/} $\HH^*(G/B)$ is free $\Z$-module linearly generated 
by the {\it Schubert classes\/} $\s_w$ labeled by the elements $w\in W$
of the Weyl group.
The {\it Poincar\'e duality\/} preserves the basis of Schubert classes:
$$
\int \s_u\cdot \s_v =\delta_{u,\wnot v}\,.\qquad
\textrm{(Kronecker's delta)}
$$
Here and everywhere below $\wnot$ is the unique Weyl group element
of maximal possible length, called the {\it longest element}.

The Schubert classes $\s_{s_i}$ generate the cohomology ring of $G/B$.
According to {\it Borel's theorem}, 
the maps $\s_{s_i}\mapsto \lambda_i$ extends to the canonical
isomorphism:
\begin{equation}
\HH^*(G/B,\Q)\simeq \mathrm{Sym}(\h_\Q^*)/I_W,
\label{eq:Borel}
\end{equation}
where $\mathrm{Sym}(\h_\Q^*)$ is the symmetric algebra of 
the $\Q$-span $\h_\Q^*\subset\h^*$ of the fundamental weights $\lambda_i$
and $I_W$ is the ideal in this algebra generated by $W$-invariant elements 
without constant term.

Let $\Z[q]=\Z[q_1,\dots,q_r]$. 
The small {\it quantum cohomology ring\/} $\QH^*(G/B)$ of $G/B$ equals, 
as a $\Z[q]$-module, to the tensor product $\HH^*(G/B)\otimes_\Z\Z[q]$.
Thus the Schubert classes $\s_w$, $w\in W$, form a $\Z[q]$-linear
basis of $\QH^*(G/B)$.
The multiplicative structure of the quantum cohomology is 
a deformation of the usual product in $\HH^*(G/B)$.
We will use the symbol ``$*$'' to denote the quantum product,
i.e., the product in the quantum cohomology ring.
The structure constants of the quantum cohomology 
are 3-point {\it Gromov-Witten invariants}:
$$
\sigma_u*\sigma_v=\sum_{w,\,d}
\<\sigma_u,\sigma_v,\sigma_{\wnot w}\>_d\,
q^d\, \sigma_w
$$
where the sum is over $w\in W$ and $d=(d_1,\dots,d_r)\in\Z_{\ge 0}^r$,
and $q^d=q_1^{d_1}\cdots q_r^{d_r}$.
The Gromov-Witten invariants
$\<\sigma_u,\sigma_v,\sigma_{\wnot w}\>_d$ are nonnegative integers
that count the
numbers of certain rational curves in $G/B$.
Their geometric definition implies that the quantum product is a commutative
and associative operation.

For a root $\alpha\in\Phi_+$
such that $h_\alpha= d_1\,h_{\alpha_1}+\cdots+d_r\,h_{\alpha_r}$, 
let $q^{h_\alpha}=q_1^{d_1}\cdots q_r^{d_r}$.
Let us assume that the variables $q_i$ are of degree 2.
Thus $\deg(q^{h_\alpha})=2\,|h_\alpha|$,
where $|h_\alpha|=(\lambda_1+\cdots+\lambda_r)(h_\alpha)$ is the 
{\it height\/} of the coroot $h_\alpha$.

Let us define, for a positive root $\alpha\in\Phi_+$, the operator 
$\T_\alpha$ acting $\Z[q]$-linearly on the quantum cohomology ring
$\QH^*(G/B)$ as
\begin{equation}
\T_\alpha:
\s_w \longmapsto 
\left\{
\begin{array}{cl}
\s_{ws_{\alpha}} & \textrm{if } \ell(ws_{\alpha}) = \ell(w)+1,\\[.1in]
q^{h_\alpha}\,\sigma_{ws_{\alpha}} & \textrm{if } \ell(ws_{\alpha}) = 
\ell(w)+1-2\,|h_\alpha|,
\end{array}
\right.
\label{eq:Talpha}
\end{equation}
The space $\QH^*(G/B)$ has the $\Z$-linear degree function 
$\deg(q^d\,\s_w)=2\,|d|+\ell(w)$.  Then the $\T_\alpha$ are homogeneous
operators on $\QH^*(G/B)$ of degree $1$.

The quantum cohomology $\QH^*(G/B)$ is generated, as
an algebra over $\Z[q]$, by the Schubert classes $\s_{s_i}$.
Thus the quantum product is uniquely determined by the following
{\it quantum Chevalley-Monk formula}:
\begin{equation}
\s_{s_i}*\s_w = \sum_{\alpha\in\Phi_+} \lambda_i(h_\alpha)\, 
\T_\alpha(\sigma_w)\,.
\label{eq:Chevalley}
\end{equation}
For root systems of type $A$, this formula was established in~\cite{FGP}.
For the case of an arbitrary type it was found by 
Peterson~\cite{Pet} (unpublished). 
A more general version of the quantum Chevalley-Monk formula, for 
$\QH^*(G/P)$, was given and rigorously proved by 
Fulton and Woodward~\cite{FW}.

\begin{figure}[ht]
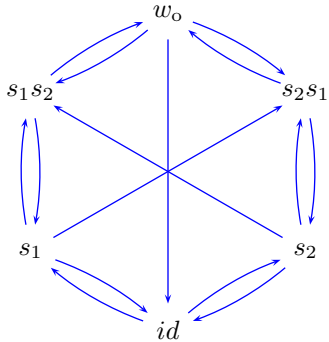

\pspicture(-40,0)(40,125)
\cnode*[linecolor=white](0,0){10}{B1}
\cnode*[linecolor=white](-52,30){10}{B2}
\cnode*[linecolor=white](52,30){10}{B3}
\cnode*[linecolor=white](-52,90){10}{B4}
\cnode*[linecolor=white](52,90){10}{B5}
\cnode*[linecolor=white](0,120){10}{B6}
\nccurve[angleA=160, angleB=-40]{->}{B1}{B2}
\nccurve[angleA=140, angleB=-20]{<-}{B1}{B2}
\nccurve[angleA=40, angleB=-160]{->}{B1}{B3}
\nccurve[angleA=20, angleB=-140]{<-}{B1}{B3}
\nccurve[angleA=100, angleB=-100]{->}{B2}{B4}
\nccurve[angleA=80, angleB=-80]{<-}{B2}{B4}
\nccurve[angleA=30, angleB=-150]{->}{B2}{B5}
\nccurve[angleA=150, angleB=-30]{->}{B3}{B4}
\nccurve[angleA=100, angleB=-100]{->}{B3}{B5}
\nccurve[angleA=80, angleB=-80]{<-}{B3}{B5}
\nccurve[angleA=40, angleB=-160]{->}{B4}{B6}
\nccurve[angleA=20, angleB=-140]{<-}{B4}{B6}
\nccurve[angleA=160, angleB=-40]{->}{B5}{B6}
\nccurve[angleA=140, angleB=-20]{<-}{B5}{B6}

\nccurve[angleA=-90, angleB=90]{->}{B6}{B1}
\rput(0,0){$id$}
\rput(52,30){$s_2$}
\rput(-52,30){$s_1$}
\rput(52,90){$s_2s_1$}
\rput(-52,90){$s_1s_2$}
\rput(0,120){$\wnot$}
\endpspicture
\caption{Quantum Bruhat graph $\Gamma_\Phi$ of type $A_2$}
\label{fig:Gamma3}
\end{figure}

Motivated by the quantum Chevalley-Monk formula~(\ref{eq:Chevalley}), 
let us define the {\it quantum Bruhat graph\/} $\Gamma_\Phi$ 
as the following directed graph on the Weyl group elements $w\in W$
with weighted edges.
Two elements $u,v\in W$ are connected by a directed edge
$e:u\to v$ if and only if $v=u s_{\alpha}$ and
one of the following two conditions
is satisfied:
$$
\ell(v)=\ell(u)+1\quad\textrm{or}\quad
\ell(v)=\ell(u)+1-2\,|h_\alpha|.
$$
If $\ell(v)=\ell(u)+1$ then weight of the edge $e$ equals 1,
and if $\ell(v)=\ell(u)+1-2\,|h_\alpha|$ 
then weight of the edge $e$ equals $q^{h_\alpha}$.
Weight of a directed path in the graph $\Gamma_\Phi$ is the product of weights
of its edges.

The graph $\Gamma_\Phi$ was investigated in~\cite{BFP} and, 
for type $A$, in~\cite{P2}.
The upward edges in $\Gamma_\Phi$, i.e., the edges that increase the length
by 1, 
are exactly the covering relations in the Bruhat order on the 
Weyl group $W$.  The downward edges correspond to additional ``quantum'' 
terms that appear in the quantum Chevalley-Monk formula~(\ref{eq:Chevalley}).  
This is why we call $\Gamma_\Phi$ the quantum Bruhat graph.

Let us say that a directed path in $\Gamma_\Phi$ from
$u$ to $v$ is {\it shortest\/} if it has the minimal possible length
among all directed paths from $u$ to $v$.

\begin{lemma}  Let $u,v\in W$ be any two Weyl group elements.
{\rm (1)} There exists a directed path from $u$ to $v$ in the graph 
$\Gamma_\Phi$.
{\rm (2)} All shortest paths from $u$ to $v$ have the same weight 
$q^{\dmin}$.
{\rm (3)} Weight of any path from $u$ to $v$ is divisible by 
$q^{\dmin}$.
\label{le:paths}
\end{lemma}

This lemma says that, for any $u,v\in W$, there is a well-defined 
{\it minimal degree\/}
$\dmin=\dmin(u,v)\in\Z_{\geq 0}^r$ such that $q^{\dmin}$ 
is the minimal possible weight of a directed path from $u$ to $v$.
The lemma easily follows from results of~\cite{BFP}.

\begin{proof}
Suppose that $a$ covers $b$ in 
the {\it weak\/} Bruhat order, i.e, $a,b\in W$, $a=b\,s_i$, 
and $\ell(a)=\ell(b)+1$,
then both directed edges $a\to b$ and $b\to a$ are present 
in the quantum Bruhat graph $\Gamma_\Phi$.  Thus any two Weyl group 
elements $u$ and $v$ can be connected by a directed path in $\Gamma_\Phi$.

According to~\cite[Lemma~6.7]{BFP}, for a
path $a\to b \to c$ of length 2 in the graph $\Gamma_\Phi$
with $a\ne c$, there exists a unique path $a\to b' \to c$ such that
$b'\ne b$.
The proof of Theorem~6.4 in~\cite{BFP}
implies that any two shortest paths from $u$ to $v$ can be obtained from 
each other by applying several switches of this type to pairs
of consecutive edges in a path. 
Moreover, any path from $u$ to $v$ can be reduced to a shortest path by applying
a sequence of such switches and cancellations of pairs of opposite 
edges $a\to b\to a$, see~\cite[Section~6]{BFP}.  

These switches of pairs of consecutive edges in a path preserve its weight.
Indeed, it is enough to verify the statement only for 
3 types $A_2$, $B_2$, and $G_2$ of rank 2 root systems, 
cf.\ proof of~\cite[Lemma~6.7]{BFP}, 
which can be easily done by a direct observation.
Thus all shortest paths from $u$ to $v$ have the same weight $q^{\dmin}$,
and weight of any path from $u$ to $v$ is divisible by $q^{\dmin}$.
\end{proof}

Let $\ell(u,v)$ be the length of a shortest path from $u$ to $v$
in the graph $\Gamma_\Phi$.  

\begin{theorem}  
Fix two Weyl group elements $u,v\in W$.  
For any $w\in W$, the coefficient of
$\s_v$ in $\s_u * \s_w$ is divisible by $q^{\dmin(u,v)}$.  
There exists $w\in W$ such that the coefficient of $\sigma_v$ in
$\sigma_u * \sigma_w$ equals $q^{\dmin(u,v)}$ times a nonzero integer.
Moreover, for any such $w$, we have $\ell(w) = \ell(u,v)$.
\label{th:minimal}
\end{theorem}

\begin{proof}
Let $S = \sigma_{s_1} + ... + \sigma_{s_r}$.  
The expression $S^{*l}=S*\cdots*S$ expands as a {\it positive\/}
integer combination of the terms $q^d\,\sigma_w$.
For any $w\in W$ of length $\ell(w) = l$, 
the expansion of $S^{*l}$ contains the class $\sigma_w$ with a strictly 
positive integer coefficient (without any quantum parameters $q_i$).
According to the 
quantum Chevalley-Monk formula~(\ref{eq:Chevalley}), the coefficient of 
$\sigma_v$ in $S^{*l}*\sigma_u$ is given by the sum over directed paths 
from $u$ to $v$, and each path comes with its weight times some 
positive integer.  By Lemma~\ref{le:paths}, 
all these terms, including the contribution of
$\sigma_w*\sigma_u$, should be divisible by $q^{\dmin(u,v)}$.
This proves the first claim of the theorem.
 
On the other hand, suppose that $l=\ell(u,v)$.  
Then there exists at least one path from $u$ to $v$ of length $l$.  
By Lemma~\ref{le:paths} 
and the quantum Chevalley-Monk formula~(\ref{eq:Chevalley}),
the coefficient of $\sigma_v$ in $S^{*l}*\sigma_u$ is equal to 
$q^{\dmin(u,v)}$ times a strictly positive integer.  The expansion of $S^{*l}$ 
in the basis of Schubert classes may also contain terms $q^d\,\sigma_{w'}$ 
with non-zero $d$.  
For all such terms, we have $\ell(w')=l'<l$ because 
$l=\deg(S^{*l})=\deg(q^d\,\sigma_{w'})=2\,|d|+l'$.  
Thus $\sigma_{w'}$ appears with a strictly positive coefficient in 
the expansion of $S^{*l'}$.  But, the coefficient of 
$\sigma_v$ in $S^{*l'}*\sigma_u$ is 0 
because there are no paths from $u$ to $v$ of length $l' < l$.  
Thus the terms $q^d\,\sigma_{w'}$ in $S^{*l}$ cannot make any contribution
to the (nonzero) coefficient of $\s_v$ in $S^{*l}*\s_u$.

This means that there exists at least one term 
$\mathrm{Const}\cdot \sigma_w$ in 
the expansion of $S^{*l}$ with $\ell(w)=\ell(u,v)$ 
that {\it do make\/} a contribution to the coefficient of 
$\s_v$ in $S^{*l}*\s_u$.
In other words, the coefficient of $\sigma_v$ in $\sigma_w * \sigma_u$ equals 
$q^{\dmin(u,v)}$ times a nonzero integer, which proves 
the second claim of the theorem.
\end{proof}
 
The next claim strengthens a result due to 
Fulton and Woodward~\cite[Theorem~9.1]{FW}.
Let us say that a monomial $q^d$ {\it occurs\/} 
in the quantum product $\sigma_u*\sigma_v$
of two Schubert classes, if there exists $w\in W$
such that the Gromov-Witten invariant $\<\s_u,\s_v,\s_w\>_d$ is nonzero.
Theorem~\ref{th:minimal} can be equivalently reformulated as follows.

\begin{corollary}
For any pair $u,v\in W$, the monomial $q^{\dmin(u,\wnot v)}$ is the 
{\it unique\/} minimal monomial that occurs in the quantum 
product $\sigma_u*\s_v$.
\end{corollary}

Note that Fulton and Woodward~\cite{FW} described the set of minimal
monomials that occur in $\sigma_u*\s_v$, but their construction does 
not immediately imply that this set consists of a single element.

It would be interesting to describe {\it all\/} monomials $q^d$
that occur in a quantum product $\s_u*\s_v$.
Theorem~\ref{th:minimal} suggests that such monomials $q^d$ 
should be weights of paths from $u$ to $\wnot v$ 
in the graph $\Gamma_\Phi$ that satisfy certain additional condition.
In the following sections we describe such a class of paths
for type $A$ root systems.

\section{Flag manifold and Schubert polynomials}

There are some notions and results related to (quantum) cohomology,
which are peculiar for type $A$.  In this section we 
discuss this type $A$ theory, which  includes the Cauchy identity 
for Schubert polynomials, see~\cite{Mac}, 
the elementary quantization rule, and 
the quantum Pieri operators, see \cite{FGP, P1} for more details.
\medskip

For type $A_{n-1}$, the homogeneous space $G/B$ is the 
complex {\it flag manifold\/}
$\Fl_n=\mathrm{SL}(n,\C)/B$.
The corresponding Weyl group is the the symmetric group $W=S_n$ 
of order $n$ permutations.
Thus the Schubert classes $\s_w$, which form a $\Z$-basis of 
the cohomology $\HH^*(\Fl_n)$, are labeled by permutations $w\in S_n$.

In this case, Borel's theorem~(\ref{eq:Borel}) implies that
the cohomology ring is isomorphic to the following quotient
of the polynomial ring:
\begin{equation}
\HH^*(\Fl_n)\simeq \Z[x_1,\dots,x_n]/\<e_1,\dots,e_n\>,
\label{eq:Borel-A}
\end{equation}
where $e_i=e_i(x_1,\dots,x_n)$ are the elementary symmetric polynomials.

Lascoux and Sch\"utzenberger~\cite{LSc},
using constructions of Bernstein-Gelfand-Gelfand and Demazure, 
defined the {\it Schubert polynomials\/} 
$\S_w(x)=\S_w(x_1,\dots,x_n)$ in the polynomial ring $\Z[x_1,\dots,x_n]$,
which are particularly nice polynomial representatives of the 
Schubert classes $\s_w$.  
These polynomials have nonnegative integer coefficients
and are stable under the standard embedding $S_n\hookrightarrow S_{n+1}$.
Another important property of the Schubert polynomials is
the following {\it Cauchy identity}, see~\cite{Mac}:
\begin{equation}
\prod_{i+j\leq n} (x_i+y_j)
=
\sum_{w\in S_n} \S_w(x_1,\dots,x_n)\cdot \S_{w\wnot}(y_1,\dots,y_n).
\label{eq:Cauchy}
\end{equation}
For type $A_{n-1}$, the longest permutation in $S_n$ is given by 
$\wnot=n\,\dots 2\,1$.

The left-hand side of~(\ref{eq:Cauchy})
is the double Schubert polynomial
$\S_\wnot(x,-y)$.  It can be expanded in terms of elementary 
symmetric polynomials as
\begin{equation}
\prod_{i+j\leq n} (x_i+y_j)=
\prod_{k=1}^{n-1}\,\sum_{i=0}^{k} y_{n-k}^{k-i}\, e_i(x_1,\dots,x_{k}).
\label{eq:rhs-cauchy}
\end{equation}

As a linear space, the quantum cohomology ring $\QH^*(\Fl_n)$ 
of the flag manifold equals $\HH^*(\Fl_n)\otimes\Z[q]$,
where $\Z[q]=\Z[q_1,\dots,q_{n-1}]$. 

For type $A_{n-1}$, the definition~(\ref{eq:Talpha}) 
of the operators $\T_\alpha$ can be written as follows.
For $1\leq i<j\leq n$, let $\T_{ij}$ be the operator
that acts $\Z[q]$-linearly on the quantum cohomology 
$\QH^*(\Fl_n)$ by the formula
\begin{equation}
\T_{ij}:\sigma_w \longmapsto 
\left\{
\begin{array}{cl}
\sigma_{ws_{ij}} & \textrm{if } \ell(ws_{ij}) = \ell(w)+1,\\[.1in]
q_{ij}\,\sigma_{ws_{ij}} & \textrm{if } \ell(ws_{ij}) = \ell(w)+1-2(j-i),
\end{array}
\right.
\label{eq:Tij}
\end{equation}
where $\ell(w)$ denote the length of permutation $w$,
$s_{ij}\in S_n$ is the transposition of $i$ and $j$, and 
$q_{ij} = q_i q_{i+1}\cdots q_{j-1}$.
The Coxeter generators of $S_n$ are the adjacent transpositions
$s_k=s_{k\,k+1}$, $k=1,\dots,n-1$.  The quantum Chevalley-Monk
formula~(\ref{eq:Chevalley})
specializes to the following formula proved in~\cite{FGP}:
\begin{equation}
\sigma_{s_k}*\sigma_w = \sum_{i\leq k<j} \T_{ij}(w),
\label{eq:Monk}
\end{equation}
for any $w\in S_n$ and $k=1,\dots,n-1$.

Let us define an involution $\omega$ on the quantum cohomology ring
$\QH^*(\Fl_n)$ by setting
$\omega: f(q_1,\dots,q_{n-1})\mapsto f(q_{n-1},\dots,q_1)$
and 
$$
\omega:q^d\sigma_w \mapsto \omega(q^d)\,\sigma_{\wnot w \wnot}\,,
$$
and extending it by linearity.  Easy observation shows that the quantum
Chevalley-Monk formula~(\ref{eq:Monk}) is invariant under $\omega$.
Thus $\omega$ is an automorphism of $\QH^*(\Fl_n)$.

Let $e_i^{(k)}=\sigma_{(k,k+1,\dots,k+i)}$ and
$h_i^{(k)}=\sigma_{(k+i,\dots,k+1,k)}$ 
be the Schubert classes represented in~(\ref{eq:Borel-A}) 
by the elementary and complete homogeneous symmetric functions 
in the first $k$ variables: 
$e_i(x_1,\dots,x_k)=\S_{(k,k+1,\dots,k+i)}$ and
$h_i(x_1,\dots,x_k)=\S_{(k+i,\dots,k+1,k)}$. 
Here we use cycle notation for permutations:
$(k,k+1,\dots,k+i)=s_k s_{k+1} \cdots s_{k+i-1}$ 
and $(k+i,\dots,k+1,k)=s_{k+i-1}\cdots s_{k+1}s_{k}$.
The involution $\omega$ on $\HH^*(\Fl_n)$
switches these two families of Schubert classes 
$e_i^{(k)}\stackrel{\omega}\longleftrightarrow h_i^{(n-k)}$.

Theorem~1.1 from \cite{FGP} is equivalent to saying that
the quantum product of the Schubert classes
$e_{i_1}^{(1)},e_{i_2}^{(2)},\dots, e_{i_{n-1}}^{(n-1)}$
in the ring $\QH^*(\Fl_n)$ is exactly the same as the 
classical product of these Schubert classes in $\HH^*(\Fl_n)$:
\begin{equation}
e_{i_1}^{(1)}*e_{i_2}^{(2)}*\cdots* e_{i_{n-1}}^{(n-1)}
=e_{i_1}^{(1)}\cdot e_{i_2}^{(2)}\cdot\, \cdots \,\cdot e_{i_{n-1}}^{(n-1)}.
\label{eq:FGP}
\end{equation}
We call this identity the {\it elementary quantization rule}.
Applying the involution $\omega$ to both sides, we get a similar
statement for the $h_i^{(k)}$: 
$$
h_{i_1}^{(1)}*h_{i_2}^{(2)}*\cdots* h_{i_{n-1}}^{(n-1)}
=h_{i_1}^{(1)}\cdot h_{i_2}^{(2)}\cdot\, \cdots \,\cdot h_{i_{n-1}}^{(n-1)}.
$$

Let us define the {\it quantum Pieri operators\/} 
$\E_i^{(k)}$ and $\H_i^{(k)}$,  acting
on the quantum cohomology $\QH^*(\Fl_n)$, as 
\begin{equation}
\begin{array}{l}
\ds
\E_i^{(k)} = \sum \T_{a_1\,b_1}\cdots T_{a_i\,b_i},\\[.1in]
\ds
\H_i^{(k)} = \sum \T_{c_1\,d_1}\cdots T_{c_i\,d_i},
\end{array}
\label{eq:Pieri-operators}
\end{equation}
where the first sum is over $a_1,\dots,a_i,b_1,\dots,b_i$
such that 
$$
a_1,\dots, a_i\leq k< b_1\leq \dots\leq b_i
\textrm{ and } a_1,\dots,a_i \textrm{ are distinct}; 
$$
and the second sum is over $c_1,\dots,c_i,d_1,\dots,d_i$ such that
$$
c_1\leq \dots\leq c_i\leq k< d_1, \dots, d_i
\textrm{ and } d_1,\dots,d_i \textrm{ are distinct}. 
$$

In~\cite{P1} we showed how to deduce from the quantum Chevalley-Monk 
formula~(\ref{eq:Monk})
the following {\it quantum Pieri formulas\/} for the quantum 
product of $e_i^{(k)}$ and $h_i^{(k)}$ with any Schubert class:
\begin{equation}
\begin{array}{l}
e_i^{(k)}*\s_w = \E_i^{(k)}(\s_w),\\[.1in]
h_i^{(k)}*\s_w = \H_i^{(k)}(\s_w).
\end{array}
\label{eq:Pieri}
\end{equation}
In a different form, the quantum Pieri formulas were earlier 
given by Ciocan-Fontanine in~\cite{CF}.
In the proof of these formulas given in~\cite{P1} we only used the 
quadratic relations for the operators $\T_{ij}$ found by Fomin
and Kirillov~\cite{FK2}.

\section{Path Schubert polynomials}
\label{sec:path=schubert}

In this section we define path Schubert polynomials and establish
their relation with the Gromov-Witten invariants for the flag manifold.
These polynomials can be expressed in terms of paths in the quantum
Bruhat graph.  As a corollary, we obtain a combinatorial rule for
all monomials $q^d$ that occurs in the quantum product $\s_u*\s_v$
in terms of paths.
\medskip

Let us define the operator $\H(y)$, acting on the space 
$\QH^*(\Fl_n)\otimes\Z[y_1,\dots,y_n]$, as the following 
combination of quantum Pieri operators:
$$
\H(y)=\sum_{\beta\leq\delta}y^{\delta-\beta} \,
\H_{\beta_1}^{(1)}\cdots \H_{\beta_{n-1}}^{(n-1)}
$$
where the sum is over nonnegative $\beta=(\beta_1,\dots,\beta_{n-1})$,
$y^{\beta} = y_1 ^{\beta_1} y_2 ^{\beta_2}\cdots y_{n-1} ^{\beta_{n-1}}$,
and $\delta=(n-1,n-2,\dots,1)$.

Let us define, for all pairs $u,v\in S_n$, 
the {\it path Schubert polynomials\/} $\S_{u,v}= \S_{u,v}(y,q)$ 
in the polynomial ring $\Z[y_1,\dots,y_n,q_1,\dots,q_{n-1}]$
as the matrix elements of the operator $\H(y)$ in the basis
of Schubert classes:
$$
\H(y):\s_u\longmapsto \sum_{v\in S_n}\S_{u,v}\cdot \sigma_v.
$$
It follows from the definitions that the $\S_{u,v}$ are polynomials
with nonnegative integer coefficients.

The specialization of the $\S_{u,v}$ for $q_1=\dots=q_{n-1}=0$ are
exactly the {\it skew Schubert polynomials\/} $\S_{v/u}\in\Z[x_1,\dots,x_n]$
recently introduced by Lenart and Sottile~\cite{LS}.
Remark that the quantum Schubert polynomials $\S_w^q$ from~\cite{FGP}
{\it cannot\/} be obtained by a specialization of 
the path Schubert polynomials $\S_{u,v}$.  Indeed, 
all coefficients of $\S_{u,v}$ are nonnegative, whereas
$\S_w^q$ may have negative coefficients.
Nevertheless, as we will see, the polynomials $\S_{u,v}$ are intimately 
related to the quantum cohomology $\QH^*(\Fl_n)$.

The definition of path Schubert polynomials $\S_{u,v}$ 
can be reformulated in terms of paths in the quantum Bruhat graph.
Let $\Gamma_n=\Gamma_\Phi$ be the quantum Bruhat graph for
type $A_{n-1}$ root system $\Phi$.
The graph $\Gamma_n$ is a directed graph on the set of permutations in $S_n$.
Each edge $e$ in $\Gamma_n$ is assigned a certain weight and a certain label.
Two permutations $u$ and $v$ are connected by an directed edge
$e:u \stackrel{i\,j}\longrightarrow v$ labeled $(i,j)$,
$1\leq i<j\leq n$, if and only if $v=u s_{ij}$ and
one of the following two conditions
is satisfied:
$$
\ell(v)=\ell(u)+1\quad\textrm{or}\quad
\ell(v)=\ell(u)+1-2(j-i).
$$
If $\ell(v)=\ell(u)+1$ then weight of the edge $e$ equals 1,
and if $\ell(v)=\ell(u)+1-2(j-i)$ then weight of the edge $e$ equals 
$q_{ij}=q_i q_{i+1}\cdots q_{j-1}$.
Weight of a directed path in the graph $\Gamma_n$ is the product of weights
of its edges.

Let us say that a directed path
$u_0 \stackrel{a_1\,b_1}\longrightarrow u_1 \stackrel{a_2\,b_2} \longrightarrow
\cdots \stackrel{a_i\,b_i}\longrightarrow u_i$ 
in the graph $\Gamma_n$
is {\it $y_k^i$-admissible\/} if 
$a_1\leq \dots\leq a_i\leq k< b_1, \dots, b_i$
and $b_1,\dots,b_i$ are distinct, cf.\ the 
definition~(\ref{eq:Pieri-operators})
of quantum Pieri operators $\H_i^{(k)}$.
More generally, let us say that a path $P$ in the graph $\Gamma_n$ is 
{\it $y^{\beta}$-admissible}, for a monomial
$y^\beta=y_1^{\beta_1}\dots y_{n-1}^{\beta_{n-1}}$, if 
the path $P$ is a concatenation of $n-1$ paths $P=P_1\circ\dots\circ P_{n-1}$
such that $P_k$ is $y_k^{\beta_k}$-admissible for $k=1,\dots,n-1$.

By definition, the path Schubert polynomial $\S_{u,v}$ is equal to
$$
\S_{u,v} = \sum_{\beta} y^{\delta-\beta}\sum_{P}  \weight(P)
$$
where the second summation is over all $y^\beta$-admissible paths $P$
in the graph $\Gamma_n$ with initial vertex $u$ and terminal vertex $v$.

\medskip
Every polynomial in the $y_i$ with $\deg_{y_i}\leq n-i$
can be expressed as a linear combination of the usual
Schubert polynomials $\S_w(y)=\S_w(y_1,\dots,y_n)$.  
The next theorem claims
that the coefficients in such an expansion for 
path Schubert polynomials are exactly the Gromov-Witten invariants.

\begin{theorem}
For any $u,v\in S_n$, we have
$$
\S_{u,\wnot v}(y,q) = \sum_{w,\,d}\<\sigma_u,\sigma_v,\sigma_{\wnot w}\>_d \,
q^d\, \S_w(y),
$$
where the sum is over $w\in S_n$ and $d\in\Z_{\geq 0}^{n-1}$.  
In other words, the coefficient of $\S_w$ in 
the Schubert-expansion of $\S_{u,\wnot v}$ is equal to 
the coefficient of the Schubert class $\sigma_w$ in 
the quantum product $\sigma_u*\sigma_v$.
\label{th:main}
\end{theorem}

This statement is a generalization to the quantum cohomology 
of a recent result by Lenart and Sottile \cite[Theorem~1]{LS}.
On the other hand, it is a flag manifold analogue of 
Theorem~6.3 from~\cite{P3}.

Theorem~\ref{th:main} implies, in particular, that
$$
\S_{w,\wnot}=\S_{1,\wnot w}=\S_w.
$$
Thus all $y^\beta$-admissible paths in $\Gamma_n$ with
the initial vertex $1$ or with the terminal vertex $\wnot$
have weight 1, i.e., they are saturated chains in the Bruhat order.

Actually, Theorem~\ref{th:main} easily follows from 
the Cauchy identity~(\ref{eq:Cauchy}),
the elementary quantization rule~(\ref{eq:FGP}), and
the quantum Pieri formula~(\ref{eq:Pieri}).

\begin{proof}[Proof of Theorem~\ref{th:main}]
Taking the image of the Cauchy identity~(\ref{eq:Cauchy}), 
with the expanded left-hand side~(\ref{eq:rhs-cauchy}),
in the ring $\HH^*(\Fl_n)\otimes\Z[y_1,\dots,y_n]$, we get
$$
\sum_\beta y^{\delta-\beta}\, 
e_{\beta_1}^{(n-1)}\cdots e_{\beta_{n-1}}^{(1)}
=\sum_{w\in S_n} \S_w(y)\cdot \s_{w\wnot}.
$$
According to the elementary quantization rule~(\ref{eq:FGP}), we 
can write the following identity in the ring 
$\QH^*(\Fl_n)\otimes\Z[y_1,\dots,y_n]$:
$$
\sum_\beta y^{\delta-\beta}\, 
e_{\beta_1}^{(n-1)}*\cdots *e_{\beta_{n-1}}^{(1)}
=\sum_{w\in S_n} \S_w(y)\cdot \s_{w\wnot}\,.
$$
Applying the involution $\omega:\QH^*(\Fl_n)\to\QH^*(\Fl_n)$, 
we get
$$
\sum_\beta y^{\delta-\beta} \, 
h_{\beta_1}^{(1)}*\cdots *h_{\beta_{n-1}}^{(n-1)}
=\sum_{w\in S_n} \S_w(y)\cdot \s_{\wnot w}\,.
$$
Since the operators of quantum multiplications by the $h_i^{(k)}$
are given by the quantum Pieri operators $\H_i^{(k)}$, the previous identity
implies that 
$$
\H(y): \s_u\mapsto \sum_{w\in S_n} \S_{w}(y)\cdot (\sigma_{\wnot w}*\sigma_u).
$$
The last formula is equivalent to the claim of theorem.
\end{proof}

Remark that the Cauchy identity was an essential
ingredient in the approach of Kirillov and Maeno~\cite{KM} to the quantum 
cohomology.

Theorem~\ref{th:main} implies a combinatorial description for the
set of monomials $q^d$ in the quantum parameters that occurs 
in the quantum product of two Schubert classes.
Let us say that a directed path 
$u_0 \stackrel{a_1\,b_1}\longrightarrow u_1 \stackrel{a_2\,b_2} 
\longrightarrow \cdots \stackrel{a_l\,b_l}\longrightarrow u_l$ 
in the graph $\Gamma_n$ is {\it admissible\/}
if there exists a sequence $k_1\leq \dots \leq k_l$ such
that $a_i\leq k_i < b_i$, for $i=1,\dots,l$, and all pairs
$(k_1,b_1),\dots,(k_l,b_l)$ are distinct.

\begin{corollary}
For any $u, v\in S_n$ and $d\in\Z_{\geq0}^{n-1}$, the monomial $q^d$ 
occurs in $\sigma_u*\sigma_v$ if and only if there exists an admissible 
path of weight $q^d$ from $u$ to $\wnot v$ in the graph~$\Gamma_n$.
\end{corollary}

\begin{proof}
Let $\S_{u,v,d}\in\Z[y_1,\dots,y_n]$ be the coefficient 
of $q^d$ in the path Schubert polynomial $\S_{u,v}$.
By Theorem~\ref{th:main}, a monomial $q^d$ occurs in $\s_u*\s_v$
if and only if $\S_{u,\wnot v,d}$ is non-zero.
On the other hand, by the definition of path Schubert polynomials, 
$\S_{u,\wnot v,d}$ is non-zero if and only if there exists a 
path from $u$ to $\wnot v$ of weight~$q^d$, which is 
$y^\beta$-admissible, for some $\beta$.  
Every $y^\beta$-admissible path is admissible and every admissible path
reduces to a $y^\beta$-admissible path.
\end{proof}

\medskip

I would like to thank Misha Kogan and Chris Woodward for 
interesting discussions and helpful correspondence.

\end{document}